\documentclass[10pt,draft,reqno]{amsart}
     \makeatletter
     \def\section{\@startsection{section}{1}%
     \z@{.7\linespacing\@plus\linespacing}{.5\linespacing}%
     {\bfseries%\normalfont\scshape
     \centering
     }}
     \def\@secnumfont{\bfseries}
     \makeatother
\newtheorem{theorem}{Theorem}[section]

\newtheorem{proposition}[theorem]{Proposition}
\newtheorem{corollary}[theorem]{Corollary}
\theoremstyle{definition}
\newtheorem{definition}[theorem]{Definition}
\newtheorem{example}[theorem]{Example}
\theoremstyle{remark}
\newtheorem{remark}[theorem]{Remark}
\numberwithin{equation}{section}
\begin{document}


%\ documentclass{amsart} 
% %\documentclass[twoside, 12pt]{article} 
% \usepackage{latexsym} 
% \usepackage{amssymb} 
% \usepackage{amsmath} 
% %\input ltex.tex 
% \font\authorfont=cmti12 
 
% \renewcommand{\baselinestretch}{1.2} 
% \renewcommand{\theequation}{\thesection.\arabic{equation}} 
% \renewcommand{\thefigure}{\thesection.\arabic{figure}} 

% \topmargin      0.25truein 
% \oddsidemargin  0.0truein 
% \evensidemargin 0.0truein 
% \textheight     8.5truein \textwidth      6.0truein
% %\textheight     8.5truein 
% %\textwidth      6.5truein 
% %\footskip       0.6truein 
% %\headheight     0.0truein 
% %\headsep        0.0truein 
% %\parskip 0.3cm 
% % \setlength{\leftmargini}{.9\leftmargini} 
 
% % \documentclass[12pt]{article} 
% % %\documentclass[twoside, 12pt]{article} 
% % \usepackage{latexsym} 
% % \usepackage{amssymb} 
% % \usepackage{amsmath} 
% % %\input ltex.tex 
% % \font\authorfont=cmti12 
 
% % \renewcommand{\baselinestretch}{1.2} 
% % \renewcommand{\theequation}{\thesection.\arabic{equation}} 
% % \renewcommand{\thefigure}{\thesection.\arabic{figure}} 
% % % 

% % %\topmargin      0.25truein 
% % %\oddsidemargin  0.0truein 
% % %\evensidemargin 0.0truein 
% % %\textheight     8.5truein 
% % %\textwidth      6.5truein 
% % %\footskip       0.6truein 
% % %\headheight     0.0truein 
% % %\headsep        0.0truein 
% % %\parskip 0.3cm 
% % % \setlength{\leftmargini}{.9\leftmargini} 

%%%%%%%%%%%%%%%%%%%%%%Definitions%%%%%%%%%%%%%%%%%%%%%%%%%%%%%%%%%%%%%%%%%%% 

%\newtheorem{theorem}{Theorem} 
\newtheorem{problem}{Problem} 
\newtheorem{conjecture}{Conjecture} 
\newtheorem{algorithm}{Algorithm} 
\newtheorem{exercise}{Exercise} 
\newtheorem{remarkk}{Remark} 
 
\newcommand{\be}{\begin{equation}} 
\newcommand{\ee}{\end{equation}} 
\newcommand{\bea}{\begin{eqnarray}} 
\newcommand{\eea}{\end{eqnarray}} 
\newcommand{\beq}[1]{\begin{equation}\label{#1}} 
\newcommand{\eeq}{\end{equation}} 
\newcommand{\beqn}[1]{\begin{eqnarray}\label{#1}} 
\newcommand{\eeqn}{\end{eqnarray}} 
\newcommand{\beaa}{\begin{eqnarray*}} 
\newcommand{\eeaa}{\end{eqnarray*}} 
\newcommand{\req}[1]{(\ref{#1})} 
 
\newcommand{\lip}{\langle} 
\newcommand{\rip}{\rangle} 

\newcommand{\uu}{\underline} 
\newcommand{\oo}{\overline} 
\newcommand{\La}{\Lambda} 
\newcommand{\la}{\lambda} 
\newcommand{\eps}{\varepsilon} 
\newcommand{\om}{\omega} 
\newcommand{\Om}{\Omega} 
\newcommand{\ga}{\gamma} 
\newcommand{\rrr}{{\Bigr)}} 
\newcommand{\qqq}{{\Bigl\|}} 
 
\newcommand{\dint}{\displaystyle\int} 
\newcommand{\dsum}{\displaystyle\sum} 
\newcommand{\dfr}{\displaystyle\frac} 
\newcommand{\bige}{\mbox{\Large\it e}} 
\newcommand{\integers}{{\Bbb Z}} 
\newcommand{\rationals}{{\Bbb Q}} 
\newcommand{\reals}{{\rm I\!R}} 
\newcommand{\realsd}{\reals^d} 
\newcommand{\realsn}{\reals^n} 
\newcommand{\NN}{{\rm I\!N}} 
\newcommand{\DD}{{\rm I\!D}} 
\newcommand{\degree}{{\scriptscriptstyle \circ }} 
\newcommand{\dfn}{\stackrel{\triangle}{=}} 
\def\complex{\mathop{\raise .45ex\hbox{${\bf\scriptstyle{|}}$} 
     \kern -0.40em {\rm \textstyle{C}}}\nolimits} 
\def\hilbert{\mathop{\raise .21ex\hbox{$\bigcirc$}}\kern -1.005em {\rm\textstyle{H}}} %Hilbert space 
\newcommand{\RAISE}{{\:\raisebox{.6ex}{$\scriptstyle{>}$}\raisebox{-.3ex} 
           {$\scriptstyle{\!\!\!\!\!<}\:$}}} % >< one above each other 
 
\newcommand{\hh}{{\:\raisebox{1.8ex}{$\scriptstyle{\degree}$}\raisebox{.0ex} 
           {$\textstyle{\!\!\!\! H}$}}} 

\newcommand{\OO}{\won} 
\newcommand{\calA}{{\mathcal A}} 
\newcommand{\calB}{{\cal B}} 
\newcommand{\calC}{{\cal C}} 
\newcommand{\calD}{{\cal D}} 
\newcommand{\calE}{{\cal E}} 
\newcommand{\calF}{{\mathcal F}} 
\newcommand{\calG}{{\cal G}} 
\newcommand{\calH}{{\cal H}} 
\newcommand{\calK}{{\cal K}} 
\newcommand{\calL}{{\cal L}} 
\newcommand{\calM}{{\cal M}} 
\newcommand{\calO}{{\cal O}} 
\newcommand{\calP}{{\cal P}} 
\newcommand{\calU}{{\mathcal U}} 
\newcommand{\calX}{{\cal X}} 
\newcommand{\calXX}{{\cal X\mbox{\raisebox{.3ex}{$\!\!\!\!\!-$}}}} 
\newcommand{\calXXX}{{\cal X\!\!\!\!\!-}} 
\newcommand{\gi}{{\raisebox{.0ex}{$\scriptscriptstyle{\cal X}$} 
\raisebox{.1ex} {$\scriptstyle{\!\!\!\!-}\:$}}} 
\newcommand{\intsim}{\int_0^1\!\!\!\!\!\!\!\!\!\sim} 
\newcommand{\intsimt}{\int_0^t\!\!\!\!\!\!\!\!\!\sim} 
\newcommand{\pp}{{\partial}} 
\newcommand{\al}{{\alpha}} 
\newcommand{\sB}{{\cal B}} 
\newcommand{\sL}{{\cal L}} 
\newcommand{\sF}{{\cal F}} 
\newcommand{\sE}{{\cal E}} 
\newcommand{\sX}{{\cal X}} 
\newcommand{\R}{{\rm I\!R}} 
\renewcommand{\L}{{\rm I\!L}} 
\newcommand{\vp}{\varphi} 
\newcommand{\N}{{\rm I\!N}} 
\def\ooo{\lip} 
\def\ccc{\rip} 
\newcommand{\ot}{\hat\otimes} 
\newcommand{\rP}{{\Bbb P}} 
\newcommand{\bfcdot}{{\mbox{\boldmath$\cdot$}}} 
 
\renewcommand{\varrho}{{\ell}} 
\newcommand{\dett}{{\textstyle{\det_2}}} 
\newcommand{\sign}{{\mbox{\rm sign}}} 
\newcommand{\TE}{{\rm TE}} 
\newcommand{\TA}{{\rm TA}} 
\newcommand{\E}{{\rm E\,}} 
\newcommand{\won}{{\mbox{\bf 1}}} 
\newcommand{\Lebn}{{\rm Leb}_n} 
\newcommand{\Prob}{{\rm Prob\,}} 
\newcommand{\sinc}{{\rm sinc\,}} 
\newcommand{\ctg}{{\rm ctg\,}} 
\newcommand{\loc}{{\rm loc}} 
\newcommand{\trace}{{\,\,\rm trace\,\,}} 
\newcommand{\Dom}{{\rm Dom}} 
\newcommand{\ifff}{\mbox{\ if and only if\ }} 
\newcommand{\nproof}{\noindent { Proof. \ }} 
\newcommand{\remarks}{\noindent {\bf Remarks:\ }} 
\newcommand{\note}{\noindent {\bf Note:\ }}

\newcommand{\boldx}{{\bf x}} 
\newcommand{\boldX}{{\bf X}} 
\newcommand{\boldy}{{\bf y}} 
\newcommand{\boldR}{{\bf R}} 
\newcommand{\uux}{\uu{x}} 
\newcommand{\uuY}{\uu{Y}} 
 
\newcommand{\limn}{\lim_{n \rightarrow \infty}} 
\newcommand{\limN}{\lim_{N \rightarrow \infty}} 
\newcommand{\limr}{\lim_{r \rightarrow \infty}} 
\newcommand{\limd}{\lim_{\delta \rightarrow \infty}} 
\newcommand{\limM}{\lim_{M \rightarrow \infty}} 
\newcommand{\limsupn}{\limsup_{n \rightarrow \infty}} 
 
\newcommand{\ra}{ \rightarrow }

\newcommand{\ARROW}[1] 
  {\begin{array}[t]{c}  \longrightarrow \\[-0.2cm] \textstyle{#1} \end{array} } 
 
\newcommand{\AR} 
 {\begin{array}[t]{c} 
  \longrightarrow \\[-0.3cm] 
  \scriptstyle {n\rightarrow \infty} 
  \end{array}} 
 
\newcommand{\pile}[2] 
  {\left( \begin{array}{c}  {#1}\\[-0.2cm] {#2} \end{array} \right) } 
 
\newcommand{\floor}[1]{\left\lfloor #1 \right\rfloor} 
 
%for doing boldface subscripts etc., e.g. $G_{\mmbox{\boldx}}$ 
\newcommand{\mmbox}[1]{\mbox{\scriptsize{#1}}} 
 
%fraction with round brackets 
\newcommand{\ffrac}[2] 
  {\left( \frac{#1}{#2} \right)} 
 
\newcommand{\one}{\frac{1}{n}\:} 
\newcommand{\half}{\frac{1}{2}\:} 
 
\def\le{\leq} 
\def\ge{\geq} 
\def\lt{<} 
\def\gt{>} 
 
%qed 
\def\squarebox#1{\hbox to #1{\hfill\vbox to #1{\vfill}}} 
\newcommand{\nqed}{\hspace*{\fill} 
           \vbox{\hrule\hbox{\vrule\squarebox{.667em}\vrule}\hrule}\bigskip} 
%%%%%%%%%%%%%%%%%%%%%%%%%%%%%%%%%%%%%%%%%%%%%%%%%%%%%%%%%%%%%%%%%%%%%%%%%%%% 
 
\title{Persistence of   invertibility 
 in the Wiener space}

\author{ A. S. \"Ust\"unel} 
%\date{ } 
\maketitle 
\noindent 
{ABSTRACT.}{\small{ Let $(W,H,\mu)$ be the classical  Wiener space, assume 
that $U=I_W+u$ is an adapted perturbation of identity where the
perturbation $u$ is an equivalence class w.r.to the Wiener
measure.  We study   several  necessary and sufficient conditions for the  almost sure
invertibility of such maps. In particular the 
subclass of these  maps who preserve the Wiener measure are
characterized  in terms of the corresponding innovation processes.  We give the following
application: let $U$ be invertible and let  $\tau$
be stopping time. Define 
$U^\tau$ as  $I_W+u^\tau$ where  $u^\tau$ is given  by
$$
u^\tau(t,w)=\int_0^t 1_{[0,\tau(w)]}(s)\dot{u}_s(w)ds\,.
$$
We prove that $U^\tau$ is also almost surely invertible. Note that  this has
immediate applications to stochastic differential equations.
}}\\ 

\vspace{0.5cm} 
\noindent

\section{Introduction} 

\noindent
This paper continues the study of the characterization of invertible
(and/or non-invertible) adapted perturbation of identity  (API for short) on the Wiener space
using the notion of the innovation process which has been developed by
Gopinath Kallianpur and his co-authors, cf.\cite{FKK}.
In \cite{ASU-2,ASU-3} we have shown some  results about the
invertibility of the adapted perturbations of identity on the
classical Wiener space. In particular, using the notion of the
innovation process associated to an API,  we have shown that the invertibility of such a
mapping is equivalent to the equality of the energy of its
perturbation to the relative entropy of the measure that it
induces. The main ingredient for all this results originates from a
result which was born from the question of representability of the
absolutely continuous probability measures as an image, or push
forward  of the Wiener
measure  under  API, which is a causal version of Monge-Kantorovitch measure transportation
theory, cf. \cite{fandu1,FUZ}. To be accurate let  $d\nu=Ld\mu$ be a probability on the
Wiener space with  $L>0$
a.s., where $\mu$ is the Wiener measure. Then there exists an API $U=I_W+u$
such that  $U\mu=\nu$ if and only if the causal estimation of $u$ w.r. to $U$ is equal to $v\circ U$,
where $v$ is the primitive of the  uniquely defined adapted process
whose Girsanov exponential is equal to $L$,
cf. \cite{ASU-2,ASU-3}. Using this result we obtain also several
interesting information about the existence of almost sure inverses of
the maps API. The list of these results is completed by proving some
more and also by proving the equivalence between them  at the beginning of
the paper, in particular we give a precise equivalence between the
existence of strong solutions of some functional stochastic
differential equations and the invertibility of the associated API,
this result finds a nice application at the end of the present paper.. 
As a byproduct of these results  we give a complete
characterization of the API's which preserve the Wiener measure in
terms of their innovation processes and show that they are closed
under composition operator. In fact, even the existence of such API's
are quite astonishing by itself in the sense that, due to the
ergodicity of the translations in the Cameron-Martin space direction,
one would expect that, at least in  the adapted case, such
transformations would be trivial, which is a totally erroneous 
intuition as we show there. We also  prove the
connections with the notion of  Girsanov measure (cf.\cite{BOOK})
associated to an API and the existence of strong solutions of
(functional) stochastic differential equations. Finally we prove that
if an API is almost surely invertible, then, its stopped version,
where the stopping occurs only in the drift, is again almost surely
invertible. Translated into the language of stochastic differential
equations, this result is rather astonishing since the stopping
operation creates a quite singular drift. Let us note that all these
results are easily extended to the infinite dimensional situations and
also to the case of the abstract Wiener spaces by using the techniques
developed in \cite{FILT} and in \cite{BOOK}.

\section{Invertibility of API's}
\noindent
Assume that  $(W,H,\mu)$ is the classical Wiener space, i.e., 
$W=C([0,1],\R^d),\,d\geq 1$, $\mu$ is the  standard Gauss measure and  $H$ 
is the  Cameron-Martin space whose  scalar product and  norm are noted
as   $(h,k)_H=\int_0^1\dot{h}_s\cdot\dot{k}_sds$
and as  $|\cdot|_H$ respectively. We note by  $(\calF_t,t\in [0,1])$
the canonical filtration of $W$ which is completed with $\mu$-null
sets. Let now  $u$ be  any $H$-valued  random variable whose Lebesgue
density is adapted (the $dt\times d\mu$-equivalence classes of such
random variables is denoted by  $L^0_a(\mu,H)$). We note by
 $\rho(\delta u)$ the  Girsanov exponential defined as
$$
\rho(\delta u)=\exp\left(\int_0^1\dot{u}_s.dW_s-\frac{1}{2}\int_0^1|\dot{u}_s|^2ds\right)
$$
In the sequel we shall denote the It\^o integral on $[0,1]$, of $\dot{u}$ with
respect to the Wiener process as $\delta u$ where $\delta$ denotes the
divergence operator defined w.r. to the Wiener measure $\mu$ (cf. for
instance  \cite{ASU,BOOK}). As an
abuse of notation we shall use again  the same notation even if $\dot{u}$ is
not square integrable w.r.to $dt\times d\mu$.

\begin{definition}
Assume that $A,B:W\to W$ are  measurable maps, we say that $A$ is a
{\em{(almost sure) right inverse}} of $B$ if 
\begin{enumerate}
\item the image of $\mu$ under $A$, denoted as $A\mu$ is absolutely
  continuous w.r.to $\mu$,
\item
$$
B\circ A(w)=w
$$
$\mu$-almost surely.
\end{enumerate}
If there is another measurable map $C:W\to W$ such that $B$ is a right
inverse to $C$ as defined above (including the absolute continuity of
$B\mu$ w.r.to $\mu$), then we say that $B$ is almost surely invertible
and in this case obviously we have $A=C$ almost surely.
\end{definition}

\begin{theorem}
\label{inverse-1}
Assume that $U=I_W+u$ is an API such that $E[\rho(-\delta
u)]=1$. Suppose that there exists a measurable map $V:W\to W$ such
that $V\circ U=I_W$ $\mu$-a.s., i.e., $V$ is a left inverse of
$U$. Then $V$ is also a right inverse, $V\mu\sim\mu$
(i.e. equivalent), it is also an API, hence of the form $V=I_W+v$ with
$v\in L^0_a(\mu,H)$. Moreover the associated stochastic processes
$(t,w)\to U(w)(t)$ and $(t,w)\to V(w)(t)$ denoted respectively  as
$(U_t(w),(t,w)\in[0,1]\times W)$ and  $(V_t(w),(t,w)\in[0,1]\times W)$
are  the
unique strong solutions  of the following stochastic differential
equations 
\begin{equation}
\label{SDEU}
dU_t=-\dot{v}_t\circ Udt+dW_t,\,\,U_0=0\,.
\end{equation}
\begin{equation}
\label{SDEV}
dV_t=-\dot{u}_t\circ Vdt+dW_t,\,\,V_0=0\,.
\end{equation}
Conversely, assume that there are  adapted process
$(U_t(w),(t,w)\in[0,1]\times W)$ and 
$(V_t(w),(t,w)\in[0,1]\times W)$ which are adapted  strong solutions  of the
equations  (\ref{SDEU}) and (\ref{SDEV}) respectively, then
$E[\rho(-\delta u)]=E[\rho(-\delta v)]=1$ and  the
corresponding API's are  almost
sure inverses  of each other.
\end{theorem}
\nproof
We have, for any $f\in C_b(W)$,  
$$
E[f\circ V]=E[f\circ V\circ U\,\rho(-\delta u)]=E[f\rho(-\delta u)]\,,
$$
hence $V\mu\sim\mu$. Let $\Om=\{w\in W:\,V\circ U(w)=w\}$, since
$\Om\subset U^{-1}(U(\Om))${\footnote{Note that $U(\Om)$ is a Souslin
    set, hence it is universally measurable}}, we have $\mu(U(\Om))=1$
and evidently $U\circ V(w)=w$ for any $w\in \Om$ and this proves the
invertibility of $U$. It is clear also that $V$ is of the form $V=I_W+v$. To show that $V$ is an API, we need to prove
that the Lebesgue density of $v$, denoted by $\dot{v}_t$ is
$\calF_t$-measurable for almost all $t\in [0,1]$. For this, note first
that $\dot{v}\circ U=\dot{u}$ $dt\times d\mu$-a.s., hence
$\dot{v}\circ U$ is adapted to the Brownian filtration, then, by
multiplying $\dot{v}$ by $1_{B_n}\circ \dot{v}$, where $B_n$ denotes
the ball of radius $n$ in $\R^d$, we may suppose that
$\dot{v}$ is bounded. Let $\eta$ be an element of $L^\infty_a(\mu,H)$,
denote by $\pi$ the operator of optional projection, then we have
\beaa
E[(\eta\circ U,v\circ U)_H\rho(-\delta u)]&=&E[(\eta,v)_H]\\
&=&E[(\eta,\pi v)_H]\\
&=&E[(\eta\circ U,(\pi v)\circ U)_H\rho(-\delta u)]
\eeaa
since $\eta$ is arbitrary, we conclude $(\pi v)\circ U=v\circ U$
$dt\times d\mu$-a.s., since $U\mu\sim\mu$, it follows that $\pi v=v$
$dt\times d\mu$-a.s. Now the processes  $(U_t,t\in[0,1])$ and
$(V_t,t\in[0,1])$ are clearly 
strong solutions  of (\ref{SDEU}) and  (\ref{SDEV}) respectively.
Conversely, any adapted  strong solutions  of
the  equations (\ref{SDEU}), (\ref{SDEV})  define API's  $U=I_W+u,V=I_W+v$ with
the property that $\rho(-\delta u)\circ V\rho(-\delta v)=\rho(-\delta
v)\circ U\rho(-\delta u)=1$, hence from the Girsanov theorem we get
$E[\rho(-\delta u)]=E[\rho(-\delta v)]=1$.
\nqed

\remark The existence of strong solutions to (\ref{SDEU}) and
(\ref{SDEV}) simultaneously implies the fact that $E[\rho(-\delta
u)]=E[\rho(-\delta v)]=1$. If we suppose that only one of them has a
strong solution, say e.g. (\ref{SDEU}) (which says that $V$ is a left
inverse) the integrability condition $E[\rho(-\delta u)]=1$ (or the
 condition  $V\mu\sim \mu$) does not
follow automatically and it  should be  added  explicitly as it is
given in the following corollary: 
\begin{corollary}
\label{inverse-cor}
Let $U=I_W+u$ be an API with $E[\rho(-\delta u)]=1$ and let  $V:W\to W$ be a
measurable map such that $V\mu\sim\mu$ and that $U\circ V=I_W$
$\mu$-a.s. Then $U$ is almost surely invertible with inverse $V$ which
is also an API and all the conclusions of Theorem \ref{inverse-1} are
also valid.
\end{corollary}
\nproof
Let $\Om=\{w:U\circ V(w)=w\}$, since $\Om\subset V^{-1}(V(\Om))$, 
$$
E[1_{V(\Om)}\circ V]=1\,,
$$
since $V\mu\sim\mu$, we have $\mu(V(\Om))=1$, hence $V\circ U=I_W$
almost surely and  Theorem \ref{inverse-1} implies that $V$ is also an
API.
\nqed

\noindent
Another version of Corollary \ref{inverse-cor}, where we do not need
to  assume
the fact that $E[\rho(-\delta u)]=1$ is given as
\begin{corollary}
\label{inverse-cor1}
Assume that $U=I_W+u,\,V=I_W+v$ are API's such that $U\circ V=I_W$
a.s. Then $E[\rho(-\delta u)]=1$ and consequently $V\circ U=I_W$ a.s.
\end{corollary}
\nproof
Since $u,v$ are both adapted, the hypothesis implies that 
$$
\rho(-\delta u)\circ V\,\rho(-\delta v)=1
$$
a.s., hence from the Girsanov theorem, we get $E[\rho(-\delta u)]=1$
and the proof follows from Corollary \ref{inverse-cor}.
\nqed

\noindent
The above  results will be used often in terms of the Lebesgue
densities of the API's under question, hence we reformulate them below
using their densities:
\begin{corollary}
\label{pratik-cor}

\begin{enumerate}
\item Assume that $U=I_W+u$ is an API such that $E[\rho(-\delta u)]=1$
  and $V=I_W+v$ with $v\in L^0(\mu,H)$ such that 
$$
\dot{u}_t+\dot{v}_t\circ U=0
$$
$dt\times d\mu$-a.s., then $V$ is also an API and it is the almost sure
inverse of $U$.
\item 
Assume that $U=I_W+u$ is an API and $V:W\to W$ a measurable map such
that $V\mu\sim\mu$ and that 
$$
\dot{v}_t+\dot{u}_t\circ V=0
$$
$dt\times d\mu$-a.s., then $V$ is also an API and it is the almost sure
inverse of $U$.
\end{enumerate}
\end{corollary}

\begin{theorem}
\label{inverse-2}
Assume that $U=I_W+u$ is an API with $E[\rho(-\delta u)]=1$, then we
have 
$$
E[\rho(-\delta u)|U]\,\frac{dU\mu}{d\mu}\circ U=1
$$
$\mu$-a.s. In particular the following equation holds true
$$
E[\rho(-\delta u)|U]=\rho(-\delta u)
$$
if and only if $U$ is a.s. invertible.
\end{theorem}
\nproof
Let us denote by $L$ the Radon-Nikodym density of $U\mu$ w.r.to
$\mu$. From the Girsanov theorem, we have 
\beaa
E[f\circ U\,L\circ U\,E[\rho(-\delta u)|U]]&=&E[f\circ U\,L\circ U\,\rho(-\delta u)]\\
&=&E[f\,L]\\
&=&E[f\circ U]
\eeaa
for any $f\in C_b(W)$, hence the first claim follows. If $U$ is almost
surely invertible, then the sigma algebra generated by $U$ is equal to
$\calF_1$, hence the equality $E[\rho(-\delta u)|U]=\rho(-\delta u)$
follows. Conversely, suppose that the latter holds, we can denote the
density $L$ as $L=\rho(-\delta v)$, with $v\in L_a^0(\mu,H)$. The
equality implies that $\dot{u}+\dot{v}\circ U=0$ $dt\times d\mu$-a.s.,
hence $V\circ U=I_W$ $\mu$-a.s., where $V=I_W+v$ and the proof follows
from Corollary \ref{pratik-cor}.
\nqed

\noindent
The following proposition whose proof follows from the Girsanov
theorem, gives a necessary and sufficient condition for a density to
be the Radon-Nikodym derivative of an API denoted by $U$ and in such a
case we say that the measure (or the density) is {\bf represented} by the
mapping $U$: 
\begin{proposition}
\label{prop-1}
Assume that $L=\rho(-\delta v)$, where $v\in L_a^0(\mu,H)$, i.e.,
$\dot{v}$ is adapted and $\int_0^1|\dot{v}_s|^2ds<\infty$ a.s.
Then there exists $U=I_W+u$, with $u:W\to H$ adapted such that
$U\mu=L\mu$ and $E[\rho(-\delta  u)]=1$ if and only if the
following condition is satisfied:
\begin{eqnarray}
\label{MA-1}
1&=&L_t\circ U \,\,E\left[\rho(-\delta u^t)|\calU_t\right]\\
&=&L_t\circ U \,\,E\left[\rho(-\delta  u)|\calU_t\right]
\end{eqnarray}
almost surely for any $t\in [0,1]$, where $u^t$ is defined as
$u^t(\tau)=\int_0^{t\wedge \tau}\dot{u}_sds$  and  $\calU_t$ is
the sigma algebra generated by $(w(\tau)+u(\tau),\,\tau\leq t)$.
\end{proposition}
 
\noindent 
Let us calculate $E[\rho(-\delta u^t)|\calU_t]=E[\rho(-\delta  u)|\calU_t]$ in terms of the
innovation process associated to $U$. Recall that the term
innovation, which originates from the filtering theory is defined
as (cf.\cite{FKK} and  \cite{BOOK})
$$
Z_t=U_t-\int_0^t E[\dot{u}_s|\calU_s]ds
$$
and it is a $\mu$-Brownian motion with respect to the filtration
$(\calU_t,t\in [0,1])$. A similar proof as the one in \cite{FKK} shows
that any martingale with respect to the filtration of $U$ can be
represented as a stochastic integral with respect to $Z$. Hence, by
the positivity assumption, $E[\rho(-\delta  u)|\calU_t]$ can be written
as an exponential martingale
$$
E[\rho(-\delta u)|\calU_t]=\exp\left(-\int_0^t(\dot{\xi}_s,dZ_s)-\frac{1}{2}\int_0^t|\dot{\xi}_s|^2ds\right)\,.
$$
Remark also that $U$ is a Wiener process under the probability
$\hat{\rho}d\mu$ where 
$$
\hat{\rho}=\exp\left(-\int_0^t(E[\dot{u}_s|\calU_s],dZ_s)-\frac{1}{2}\int_0^t|E[\dot{u}_s|\calU_s]|^2ds\right)\,,
$$
hence a  double utilization of the Girsanov theorem gives the following
explicit result:
\begin{proposition}
\label{prop-2}
We have 
\begin{equation}
\label{cond-exp} E[\rho(-\delta
u)|\calU]=\exp\left(-\int_0^1(E[\dot{u}_s|\calU_s],dZ_s)-\frac{1}{2}\int_0^1|E[\dot{u}_s|\calU_s]|^2ds\right)\,,
\end{equation}
and
\begin{equation}
\label{cond-exp1} E[\rho(-\delta
u)|\calU_t]=\exp\left(-\int_0^t(E[\dot{u}_s|\calU_s],dZ_s)-\frac{1}{2}\int_0^t|E[\dot{u}_s|\calU_s]|^2ds\right)\,,
\end{equation}
almost surely.
\end{proposition}
Combining Propositions \ref{prop-1} and \ref{prop-2}, we obtain
\begin{theorem}
\label{thm-2}
A  necessary and sufficient condition for a density $L$, represented as
$L=\rho(-\delta v)$, where $v\in L^0_a(\mu,H)$ to be  the Radon-Nikodym
density of the image of the Wiener  measure $\mu$ under  some API,
noted as  $U=I_W+u$, is that
$$
E[\dot{u}_t|\calU_t]=-\dot{v}_t\circ U
$$
$dt\times d\mu$-almost surely.
\end{theorem}

\noindent 
Now we state and prove a  main theorem (cf. also \cite{FUZ} for
related problems):
\begin{theorem}
 \label{entropy-thm}
 Assume that $u\in L^2(\mu,H)\cap L^0_a(\mu,H)$ with $E[\rho(-\delta 
 u)]=1$. Define $L$ as
$$
L=\frac{dU\mu}{d\mu}=\rho(-\delta  v)
$$
where  $v\in L_a^0(\mu,H)$ is given by the It\^o representation
theorem. The map  $U=I_W+u$ is then almost surely
invertible with its inverse  $V=I_W+v$ if and only if
$$
E[L\log L]=\frac{1}{2}E[|u|_H^2]\,.
$$
In other words, $U$ is invertible if and only if
$$
H(U\mu|\mu)=\half \|u\|^2_{L^2(\mu,H)}\,,
$$
where $H(U\mu|\mu)$ denotes the entropy of $U\mu$ with respect to
$\mu$.
\end{theorem}
\nproof
 Since $U$ represents $Ld\mu$, we have, from Theorem \ref{thm-1}, 
$E[\dot{u}_s|\calU_s]+\dot{v}_s\circ U=0$ $ds\times d\mu$-almost
surely. Hence, from the Jensen inequality $E[|v\circ U|_H^2]\leq
E[|u|_H^2]$. Moreover the Girsanov theorem gives
\beaa
2E[L\log L]&=&E[|v|_H^2 L]\\
&=&E\left[|v\circ U|_H^2\right]\\
&=&E\left[\int_0^1|E[\dot{u}_s|\calU_s]|^2ds\right]\,.
\eeaa
Hence the hypothesis implies that
$$
E[|u|_H^2]=E[\int_0^1|E[\dot{u}_s|\calU_s]|^2ds]\,.
$$
From which we deduce that $\dot{u}_s=E[\dot{u}_s|\calU_s]$
$ds\times d\mu$-almost surely. Finally we get
$\dot{u}_s+\dot{v}_s\circ U=0$ $ds\times d\mu$, which is a
necessary and sufficient condition for the invertibility of $U$ from
Corollary \ref{pratik-cor} (cf. also \cite{ASU-2}) . The necessity is
obvious.

 \nqed
\begin{corollary}
With the notations of Theorem \ref{entropy-thm}, $U$ is not invertible if and only if we
have 
$$
\half E[|u|_H^2]>H(U\mu|\mu)\,.
$$
\end{corollary} 
\remark
This result gives an enlightenment about the celebrated counter
example of Tsirelson, cf. \cite{I-W}.

\begin{corollary}
\label{cor-3}
Assume that $(U^n=I_W+u^n,\,n\geq 1)$ is a sequence of adapted and
almost surely invertible 
perturbations of identity such that $E[\rho(-\delta u^n)]=1$ for any
$n\geq 1$.  Suppose that $(U^n,\,n\geq 1)$ converges to some
$U=I_W+u$ in $L^0(\mu,W)$ with $u\in L^0_a(\mu,H)$ with  $E[\rho(-\delta u)]=1$. If 
$$
\lim_nH(U^n\mu|\mu)=H(U\mu|\mu)\,,
$$
 then $U$ is almost surely invertible. 
\end{corollary}
\nproof
From the hypothesis, it follows that $U\mu\sim \mu$. Let 
$$
L=\frac{dU\mu}{d\mu}=\rho(-\delta v)
$$
where $v\in L^0_a(\mu,H)$ is uniquely defined from the It\^o
representation theorem. We have, from Theorem \ref{thm-2}, 
$$
E[\dot{u}_t|\calU_t]+\dot{v}_t\circ U=0
$$
almost surely. Moreover, from the Fatou lemma and from the lower semi
continuity of the Cameron-Martin norm with respect to the Banach norm
of $W$, 
\beaa
\half E[|u|_H^2]&\leq&\half E[\lim\inf_n|u_n|_H^2]\\
&\leq&\half \lim\inf_nE[|u_n|_H^2]\\
&=&H(U\mu|\mu)=\half E[|\hat{u}|_H^2]
\eeaa
where $\hat{u}(t)=\int_0^tE[\dot{u}_s|\calU_s]ds$. Since $\hat{u}$ is
an orthogonal projection of $u$, it follows that $\hat{u}=u$ a.s.,
hence, due to Theorem \ref{entropy-thm}, $U$ is almost surely invertible with inverse $V=I_W+v$.
\nqed

\begin{theorem}
\label{1-more}
Assume that $(U_n=I_W+u_n,\,n\geq 1)$ is a sequence of a.s. invertible sequence of identities with $L_n=dU_n\mu/d\mu$ satisfying the following properties:
\begin{enumerate}
\item $\lim_nL_n=L$ weakly in $L^1(\mu)$.
\item There exists a measurable map $U:W\to W$ such that
$$
\lim_nE[f\circ U_n]=E[f\circ U]
$$
for any $f\in C_b(W)$.
\item
$$
\lim_nE[L_n\log L_n]=E[L\log L]\,.
$$
\end{enumerate}
Then $\frac{dU\mu}{d\mu}=L$ and $U$ is a.s. invertible.
\end{theorem}
\nproof
By writing $u=U-I_W$, from the lower semicontinuity of the Cameron-Martin norm on $W$, we see that $u$ is an $H$-valued map, besides, it follows from the hypothesis that it is adapted to the canonical filtration. Evidently 
$$
\frac{dU\mu}{d\mu}=L\,.
$$
Again from the lower semicontinuity and from the last hypothesis
\beaa
E[L\log L]&=&\lim\inf_E[L_n\log L_n]\\
&=&\lim\inf\frac{1}{2}E[|u_n|_H^2]\\
&\geq&\frac{1}{2}E[|u|_H^2]\\
&\geq&E[L\log L]
\eeaa
hence we get 
$$
E[L\log L]=\frac{1}{2}E[|u|_H^2]
$$
which is a sufficient condition for the invertibility of $U$.
\nqed

\begin{theorem}
\label{gen-density}
Let  $U=I_W+u$, $u\in L^2_a(\mu,H)$ and let us denote by $L$ the
Radon-Nikodym derivative $dU\mu/d\mu$. Assume that 
$$
H(U\mu|\mu)=\half\|u\|^2_{L^2(\mu,H)}
$$
and that
$$
E[L\log L]+E[-\log L]<\infty\,.
$$
Then $U$ is almost surely invertible.
\end{theorem}
\nproof
Since $-\log L$ is integrable, $L$ is a.s. strictly positive, hence it
can be represented as $L=\rho(-\delta v)$, where $v\in
L_a^0(\mu,H)$. We have $L\circ U\,E[\rho(-\delta u)|\calU]\leq 1$ from the
Girsanov theorem. Using the Jensen and above  inequalities  we get
\beaa
E[L\log L]&=&E[\log L\circ U]\\
&\leq&-E[\log E[\rho(-\delta u)|\calU]]\\
&\leq & E[\delta u+\half|u|_H^2]\\
&=&E[L\log L]
\eeaa
therefore 
$$
-E[\log E[\rho(-\delta u)|\calU]]=E[\log L\circ U]
$$
and this relation implies that
$$
L\circ U E[\rho(-\delta u)|\calU]=1\,,
$$
hence $E[\rho(-\delta u)]=1$ and the proof follows from Theorem \ref{entropy-thm}.
\nqed

\noindent
We shall give below another application of Theorem \ref{thm-1} which
is about the  measure preserving adapted perturbations of identity:
\begin{theorem}
\label{thm-3}
Assume that $a\in L^2_a(\mu,H)$ with $E[\rho(-\delta a)]=1$. Define
$A=I_W+a$, then $A$ preserves the Wiener measure, i.e.,  $A\mu=\mu$,
if and only if we have 
$$
E[\dot{a}_t|\calA_t]=0
$$
$dt\times d\mu$-almost surely, where $(\calA_t,t\in [0,1])$ denotes
the filtration of $A$. In particular $A$ is equal to its innovation process.
\end{theorem}
\nproof
 From the Girsanov theorem, for any $f\in C_b(W)$, we have 
\beaa
E[f\circ A]&=&E[f]\\
&=&[f\circ A\,\rho(-\delta a)]\\
&=&[f\circ A\,E[\rho(-\delta a)|\calA]]\,.
\eeaa
Hence
$$
E[\rho(-\delta a)|\calA]=1\,,
$$
and  from Theorem \ref{thm-2}, $E[\dot{a}_t|\calA_t]=0$ a.s. Letting
$Z=(Z_t)$ be the innovation process associated to $A$, we get 
\beaa
Z_t&=&A_t-\int_0^tE[\dot{a}_s|\calA_s]ds\\
&=&A_t\,.
\eeaa
\nqed

\noindent
One can construct measure preserving API's as explained  in the
following example:
\begin{example}
Let $U=I_W+u$, where $\dot{u}$ is the shift given by B. Tsirelson,
c.f. \cite{Tsi} or \cite{I-W}, then, as it is well-known, the API  $U$ is not
invertible. On the other hand, since
$u$ is bounded, $U\mu\sim \mu$. Let
$$
L=\frac{dU\mu}{d\mu}=\rho(-\delta v)
$$
where $v\in L^0_a(\mu,H)$ is uniquely defined due to the It\^o
representation theorem. Define $V$ as  $V=I_W+v$ and let 
$A=V\circ U$. From the Girsanov theorem, we have
$$
E[f\circ A]=E[f\circ V\circ U]=E[f\circ V\rho(-\delta v)]=E[f]\,,
$$
for any $f\in C_b(W)$ and $A=I_W+a=I_W+u+v\circ U$ with $\dot{a}$
adapted. Note that this subset of  API's  is
closed with respect to the composition operation.
\end{example}

\noindent
Let us recall the notion of Girsanov measure which has been  already
described in \cite{BOOK}:
\begin{definition}
Let $(\Om,\calF,\rho)$ be a probability  space on which is  given a measurable
map $T:\Om\to\Om$. A  measure $\nu$ on $(\Om,\calF)$ is called a Girsanov measure for
(or associated to) $(\rho,T)$ if $T\nu=\rho$; in other words, if 
$$
\int_\Om f\circ T\,d\nu=\int_\Om f\,d\rho
$$
for any measurable, positive function on $\Om$.
\end{definition}

\begin{theorem}
\label{Gir-measure}
Assume that $u\in L^0_a(\mu,H)$ such that $E[\rho(-\delta
u)]=1$. Let $U=I_W+u$, then  there exists a unique absolutely
continuous (w.r.to $\mu$) Girsanov measure associated to
$(U,\mu)$ if and only if $U$ is almost surely invertible.
\end{theorem}
\nproof
To show the necessity note that $\rho(-\delta u)d\mu$ and
$E[\rho(-\delta u)|\calU]d\mu$ are two absolutely continuous Girsanov
measures. The uniqueness implies that $\rho(-\delta u)=E[\rho(-\delta
u)|\calU]$ almost surely. It follows from Theorem \ref{thm-1} that
$\dot{u}_t+\dot{v}_t\circ U=0$ $dt\times d\mu$-a.s., where $\dot{v}$
is defined as $\rho(-\delta v)=dU\mu/d\mu$, hence $U$ and $V=I_W+v$
are a.s. inverse to each other from Theorem \ref{inverse-1}. To show the sufficiency, let
$d\nu=\Gamma d\mu$ be any Girsanov measure for $(\mu,U)$ where
$\Gamma\in L^1_+(\mu)$. Then, by the a.s. invertibility of $U$ we have 
$\Gamma=\rho(-\delta u)$ almost surely.
\nqed

\noindent
The next theorem summarizes the most notable results of this section
about the invertibility of the API's:
\begin{theorem}
\label{thm-1}
Assume that  $E[\rho(-\delta u)]=1$  and denote by  $U$ the mapping 
$I_W+u$. The following properties are then equivalent
\begin{enumerate}
\item  $U$ is almost surely invertible and its inverse $V$ is of the form
  $V=I_W+v$ with  $v\in L^0_a(\mu,H)$,
\item The following stochastic differential equation
\beaa
dV_t&=&-\dot{u}_t\circ V\,dt+dW_t\\
V_0&=&0
\eeaa
has a unique strong solution.
\item The following relation holds true
$$
\half \int_W|u|_H^2d\mu=\int_W \frac{dU\mu}{d\mu}\log\frac{dU\mu}{d\mu}d\mu\,.
$$
\item We have the following identity
$$
\frac{dU\mu}{d\mu}\circ U\,\rho(-\delta u)=1
$$
almost surely.
\item $U$ has a unique absolutely continuous  Girsanov measure.
\end{enumerate}
\end{theorem}
\nproof
The equivalence between (1) and (2) follows from Theorem
\ref{inverse-1} and  the one  between (1) and (3) is proved  in
Theorem \ref{entropy-thm}.  The equivalence between (1) and (4) is
given by Theorem \ref{inverse-2} and finally the one  between
(1) and (5) is given by Theorem \ref{Gir-measure}.
\nqed

\section{Invertibility  is preserved under stopping of the adapted shifts } 
\noindent
Assume that $U=I_W+u$ is an invertible adapted perturbation of
identity, whose inverse is given by $V=I_W+v$. If $a\in[0,1]$ is any
fixed number, define $U^a$ as $U^a=I_W+u^a$, where $u^a$ is defined as 
$$
uâ(t)=\int_0^t 1_{[0,a]}(s)\dot{u}_sds\,.
$$
Then it is easy to see that $U^a$ is invertible and its inverse is
given explicitly as $V^a=I_W+v^a$, where $v^a$ is defined as $u^a$
above. A natural question is whether this property persists if  we
replace the constant $a$ with a stopping time $\tau$. The next theorem
answers this question positively:
\begin{theorem}
\label{thm-4}
Let $u,v\in L^0_a(\mu,H)$ s.t. $E[\rho(-\delta u)]=1$, define
$U=I_W+u,\,V=I_W+v$. Assume that $U$ and $V$  are  a.s. inverse to each other, i.e., 
$$
U\circ V=V\circ U=I_W
$$
$\mu$-almost surely. Let $\tau$ be a stopping time w.r.to the
filtration of the canonical Brownian motion, with values in
$[0,1]$. Define $U^\tau$ as $I_W+u^\tau$, where 
$$
u^\tau(t)=\int_0^t 1_{[0,\tau]}(s)\dot{u}_sds\,.
$$
Then $U^\tau$ has also a both sided inverse $S$ of the form
$S=I_W+\alpha$, where $\alpha\in L_a^0(\mu,H)$ satisfies the following
identity:
$$
\dot{\alpha}_t=\dot{v}_t\,1_{[0,\tau\circ S]}(t)
$$
$dt\times d\mu$-almost surely. In particular we have
\beaa
\frac{dU^\tau\mu}{d\mu}&=&E\left[\frac{dU\mu}{d\mu}|\calF_{\tau\circ S}\right]\\
&=&E\left[\rho(-\delta v)|\calF_{\tau\circ S}\right]\,.
\eeaa
\end{theorem}
\nproof
It suffices to prove the existence of some $\alpha\in L^0_a(\mu,H)$
such that 
$$
\dot{u}^\tau_t+\dot{\alpha}_t\circ U^\tau=0
$$
$dt\times d\mu$-a.s. From the hypothesis $U^\tau\mu$ is equivalent to
$\mu$, let $L^\tau$ be the corresponding Radon-Nikodym density. From
the It\^o representation theorem, there exists some $\alpha\in
L^0_a(\mu,H)$ such that 
$$
L^\tau=\rho(-\delta \alpha)\,.
$$
From \cite{ASU-2}, we have 
$$
E[\dot{u}^\tau_t|\calU^\tau_t]+\dot{\alpha}_t\circ U^\tau=0
$$
a.s., where $(\calU^\tau_t,t\in [0,1])$ is the filtration of the map
$U^\tau$, i.e., 
$$
\calU^\tau_t=\sigma(U^\tau(s),s\leq t)
$$
and 
$$
U^\tau(s)(w)=W_s(w)+\int_0^s1_{[0,\tau(w)]}(r)\dot{u}_r(w)dr\,.
$$
We claim that 
\begin{equation}
\label{eqn-1}
E[\dot{u}^\tau_t|\calU^\tau_t]=E[\dot{u}^\tau_t|\calU_t]
\end{equation}
a.s., where $(\calU_t,t\in [0,1])$ is the filtration of $U$. To prove the
relation (\ref{eqn-1}), let $A$ be in $L^\infty(\mu)$, then
\beaa
E\left[A\,E[\dot{u}^\tau_t|\calU^\tau_t]\right]&=&E\left[E[A|\calU^\tau_t]\dot{u}^\tau_t\right]\\
&=&E\left[E[A|\calU^\tau_t]\,1_{\{\tau>t\}}\dot{u}_t\right]\\
&=&E\left[E[A|\calU_t]\,1_{\{\tau>t\}}\dot{u}_t\right]\\
&=&E\left[A\,E[\dot{u}^\tau_t|\calU_t]\right]\,,
\eeaa
since $A$ is arbitrary this proves the relation (\ref{eqn-1}). Since
$u$ is adapted and $U$ is invertible $\calU_t=\calF_t$, where
$(\calF_t)$ is the filtration of the Wiener process, we get 
\begin{equation}
\label{eqn-inv}
\dot{u}^\tau_t+\dot{\alpha}_t\circ U^\tau=0
\end{equation}
$dt\times d\mu$-a.s., which implies that $I_W+\alpha$ is the two-sided
inverse of $U^\tau$ almost surely. To complete the proof it suffices
to verify that $\alpha$ given as in the claim satisfies the relation
(\ref{eqn-inv}):
\beaa
\left(\dot{v}_t\,1_{[0,\tau\circ S]}(t)\right)\circ
U^\tau&=&\dot{v}_t\circ U^\tau\,1_{[0,\tau\circ S\circ U^\tau]}(t)\\
&=&\dot{v}_t\circ U^\tau\,1_{[0,\tau]}(t)\\
&=&\dot{v}_t\circ U\,1_{[0,\tau]}(t)\\
&=&-\dot{u}_t\,1_{[0,\tau]}(t)\\
&=&-\dot{u}_t^\tau
\eeaa
and this completes the proof.

\nqed

\noindent
\begin{example}
 A typical and elementary example is obtained if we
take $u$ to be constant vector field $h\in H$ and if $\tau$ is any
stopping time: then the mapping $T=I_W+h^\tau$, where
$h^\tau(t,w)=\int_0^{t\wedge \tau(w)}\dot{h}_sds$, is almost surely invertible.
\end{example}
\noindent
We can interpret this result in the language of the stochastic
differential equations (SDE)  as
\begin{corollary}
Assume that $u\in L^2_a(\mu,H)$ with $E[\rho(-\delta u)]=1$. If the
SDE 
\beaa
dV_t&=&-\dot{u}_t\circ Vdt+dW_t\\
V_0&=&0
\eeaa
has a unique strong solution, then so does also the following SDE
\beaa
dS_t&=&-(\dot{u}_t1_{[0,\tau]})\circ Sdt+dW_t\\
V_0&=&0
\eeaa
for any stopping time $\tau$. In particular we have 
$$
\frac{dU^\tau\mu}{d\mu}=E\left[\frac{dU\mu}{d\mu}|\calF_{\tau\circ S}\right]\,.
$$
\end{corollary}

\noindent
{\bf{Acknowledgement:}} Some parts of this work has been done while the
author was visiting the Departement of Mathematics of Bilkent
University, Ankara, Turkey.

\vspace{0.5cm}
\footnotesize{

A. S. \"UST\"UNEL: INSTITUT  TELECOM, TELECOM- PARISTECH, LTCI CNRS D\'EPT. INFRES, 
46, RUE BARRAULT, 75013, PARIS, FRANCE.

 {\em{E-mail address:}} ustunel@enst.fr}

\end{document}